\documentclass[12pt]{article}
\usepackage{amssymb,amsmath,mathrsfs,amsthm}
\usepackage{hyperref}
\usepackage{graphicx}
\usepackage{color}
\usepackage[OT2,OT1]{fontenc}
\usepackage{upquote}
\usepackage{setspace}
\usepackage{lipsum}
\usepackage[numbers,sort&compress]{natbib}

\usepackage{float} 

\usepackage{xcolor}
\usepackage{framed}
\definecolor{shadecolor}{cmyk}{0,0,0,0.05} 

\makeatletter
\def\blfootnote{\xdef\@thefnmark{}\@footnotetext}
\makeatother

\begin{document}

\title{\Large\bf A new form of the Machin-like formula for $\pi$ by iteration with increasing integers}

\date{17 April 2022}

\author{
{\small Sanjar M. Abrarov, Rehan Siddiqui, Rajinder K. Jagpal,} \\
{\small and Brendan M. Quine}
}

\maketitle

\begin{abstract}
We present a new form of the Machin-like formula for $\pi$ that can be generated by using iteration. This form of the Machin-like formula may be promising for computation of the constant $\pi$ due to rapidly increasing integers at each step of the iteration. The computational test we performed shows that, with an integer $k \ge 17$, the Lehmer measure remains small and practically does not increase after $18$ steps of iteration.
\vspace{0.25cm}
\end{abstract}

\blfootnote{Mathematica programs can be accessed here: \href{https://cs.uwaterloo.ca/journals/JIS/VOL25/Abrarov/supplement.txt}{\textcolor{blue}{\it supplement.txt}}}

\section{Introduction}

A remarkable discovery, made in 1706 by English astronomer and mathematician John Machin~\cite{Beckmann1971,Berggren2004,Borwein2008}
\begin{equation}\label{eq_1}
\frac{\pi}{4} = 4\arctan\frac{1}{5} - \arctan\frac{1}{239},
\end{equation}
had a great impact on the mathematical society of that time. Specifically, due to relatively rapid convergence, he was the first to calculate $100$ digits of $\pi$. Nowadays equation~\eqref{eq_1} is named as the Machin formula for $\pi$,
in his honor.

The equations of the kind~\cite{Abeles1993}
\begin{equation}\label{eq_2}
\frac{\pi }{4} = \sum\limits_{j = 1}^J\alpha_j\arctan\frac{1}{\beta_j}, \qquad\alpha_j,\beta_j \in \mathbb{R},
\end{equation}
are known to be the Machin-like formulas for $\pi$. Historically, some of the earliest formulas are
\begin{equation}\label{eq_3}
\frac{\pi}{4} = \arctan\frac{1}{2} + \arctan\frac{1}{3},
\end{equation}
\begin{equation}\label{eq_4}
\frac{\pi}{4} = 2\arctan\frac{1}{2} - \arctan\frac{1}{7},
\end{equation}
\begin{equation}\label{eq_5}
\frac{\pi}{4} = 2\arctan\frac{1}{3} + \arctan\frac{1}{7},
\end{equation}
due to Euler, Hermann, and Hutton, respectively.

Since the Maclaurin series expansion of the arctangent function is given by
\begin{equation}\label{eq_6}
\arctan x = x - \frac{1}{3}x^3 + \frac{1}{5}x^5 - \frac{1}{7}x^7 + \cdots = x + O(x^3),
\end{equation}
we can see that it is very desirable to get the arguments of the arctangent function as small as possible to improve the convergence rate in the Machin-like formula~\eqref{eq_2} for $\pi$.

In 1938 Lehmer introduced a measure~\cite{Lehmer1938,Tweddle1991,Abeles1993,Wetherfield1996}
\begin{equation}\label{eq_7}
\mu  = \sum\limits_{j = 1}^J\frac{1}{\log_{10}\left|\beta _j\right|}.
\end{equation}
showing how much computational labor is needed to compute a given Machin-like formula for $\pi$. In particular, a given Machin-like formula for $\pi$ is more efficient if its constants $\beta_j$ are larger by absolute value and if the number $J$ of the terms in equation~\eqref{eq_2} is smaller. More detailed description and significance of Lehmer's measure in the computation of constant $\pi$ by using the Machin-like formulas can be found in~\cite{Wetherfield1996}.

Application of the Machin-like formulas with small Lehmer measure is one of the most efficient ways to compute digits of $\pi$. In particular, in 2002, Kanada computed over one trillion digits of $\pi$ by using a pair of self-checking Machin-like formulas~\cite{Calcut2009,Agarwal2013}.

In 1997, Chien-Lih published a remarkable six-term Machin-like formula for $\pi$~\cite{Chien-Lih1997}
\[
\begin{aligned}
\frac{\pi }{4} &= 183\arctan\frac{1}{239} + 32\arctan\frac{1}{1023} - 68\arctan\frac{1}{5832} \\
&+ 12\arctan\frac{1}{110443} - 12\arctan\frac{1}{4841182} - 100\arctan\frac{1}{6826318}
\end{aligned}
\]
with relatively small Lehmer measure $\mu \approx 1.51244$. Later Chien-Lih showed how to reduce Lehmer's measure even further by generating the two-term Machin-like formulas for $\pi$ by iteration involving Euler-type identities. However, his method is not simple and requires some algorithmic manipulations at each step of the iteration~\cite{Chien-Lih2004}.

Previously, Abrarov et al.~\cite{Abrarov2017} developed a new method of generating the two-term Machin-like formulas for $\pi$
\begin{equation}\label{eq_8}
\frac{\pi}{4} = 2^{k - 1}\arctan\frac{1}{\beta_1} + \arctan\frac{1}{\beta_2},
\end{equation}
where the constant $\beta_1$ can be chosen as a positive integer
\begin{equation}\label{eq_9}
\beta_1 = \left\lfloor\frac{a_k}{\sqrt{2 - a_{k - 1}}}\right\rfloor
\end{equation}
such that the nested radicals are defined as $a_k = \sqrt{2 + a_{k - 1}}$ and $a_0 = 0$. Unlike Chien-Lih's method of generating the two-term Machin-like formulas for $\pi$, this method is significantly easier in implementation as it is based on simple two-step iteration~\eqref{eq_12} below. Therefore, this method is more efficient to reduce the Lehmer measure than the method proposed by Chien-Lih in his work~\cite{Chien-Lih2004}.

In general, the constants $\beta_j$ in equation~\eqref{eq_2} may not necessarily be integers. In particular, Chien-Lih~\cite{Chien-Lih2004} and Abrarov et al.~\cite{Abrarov2017} methods of generating the two-terms Machin-like formulas for $\pi$ result in quotients with rapidly increasing number of digits in their numerators and denominators. As a consequence, this makes the computation of the arctangent function difficult due to exponentiation with increasing orders. In this work, we propose a method showing how this problem can be effectively resolved by generating a new form of the Machin-like formula for $\pi$. This method is based on simple iteration and, due to rapidly increasing values of integers, may be promising for computation of $\pi$ with rapid convergence.

\section{Results and discussion}\label{RD}

\subsection{Methodology}

The first constant ${\beta _1}$ from equation~\eqref{eq_8} that was derived from the following identity~\cite{Abrarov2017}
$$
\frac{\pi}{4} = 2^{k - 1}\arctan\frac{\sqrt{2-a_{k-1}}}{a_k}, \qquad k \ge 1, \qquad k \in \mathbb{Z}.
$$
can be chosen arbitrarily. However, it would be reasonable to choose it in such a way that the ratio 
$$
\frac{2^{k - 1}}{\beta _1} \approx \frac{\pi}{4}.
$$
Therefore, if we imply that $\beta_1$ is an integer, the best choice is given either by equation~\eqref{eq_9} or by equation
\begin{equation}\label{eq_10}
\beta_1 = \left\lceil\frac{a_k}{\sqrt{2 - a_{k - 1}}}\right\rceil,
\end{equation}
based on floor and ceiling functions, respectively. Once the constant ${\beta _1}$ is chosen, the second constant can be found as~\cite{Abrarov2017}
\begin{equation}\label{eq_11}
\beta_2 = \frac{2}{\left((\beta_1 + i)/(\beta_1 - i)\right)^{2^{k - 1}} - i} - i.
\end{equation}

It is interesting to note that all four earliest two-term Machin-like formulas~\eqref{eq_1},~\eqref{eq_3},~\eqref{eq_4} and~\eqref{eq_5} for $\pi$ can be readily found with help of equations~\eqref{eq_8} and~\eqref{eq_11}. In particular, the original Machin formula~\eqref{eq_1} for $\pi$ can be obtained by substituting $k = 3$ into equations~\eqref{eq_9},~\eqref{eq_11} and then~\eqref{eq_8}. The Euler equation~\eqref{eq_3} can be obtained by substituting $k = 1$ and ${\beta_1} = 2$  into equation~\eqref{eq_11} and then~\eqref{eq_8}. The Hutton equation~\eqref{eq_4} can be found by substituting $k = 2$ into equations~\eqref{eq_9},~\eqref{eq_11} and then~\eqref{eq_8}. Finally, the Hermann equation~\eqref{eq_5} can be found by substituting $k = 2$ into equations~\eqref{eq_10},~\eqref{eq_11} and then~\eqref{eq_8}.

While equation~\eqref{eq_11} is useful when integer $k$ is small, its application becomes problematic when $k$ increases. Such a problem arises as a result of rapidly increasing power ${2^{k - 1}}$ in the denominator of equation~\eqref{eq_11}. However, as we have shown in publication~\cite{Abrarov2017}, this problem can be resolved by using simple two-step iteration
\begin{equation}\label{eq_12}
\left.
\begin{aligned}
\sigma_n &= \sigma_{n - 1}^2 - \tau_{n - 1}^2 \\
\tau_n &= 2\sigma_{n - 1}\tau_{n - 1}
\end{aligned}
\right\},\qquad n = \left\{2,3,4, \ldots, k \right\},
\end{equation}
with initial values defined as
\[
\sigma_1 = \frac{\beta_1^2 - 1}{\beta_1^2 + 1}
\]
and
\[
\tau_1 = \frac{2\beta_1}{\beta_1^2 + 1}
\]
such that the second constant can be found from the ratio
\begin{equation}\label{eq_13}
\beta_2 = \frac{\sigma_k}{1 - \tau_k}.
\end{equation}

There is an identity for the arctangent function
$$
\arctan x + \arctan y = \arctan \frac{x + y}{1 - xy}.
$$
It is not difficult to see that from this identity we get
\begin{equation}\label{eq_14}
\arctan(x + y) = \arctan x + \arctan\frac{y}{1 + (x + y)x}.
\end{equation}
	
Let $z \in \mathbb{R}$ and $z = x + y$ such that its integer and fractional parts are given by $x = \left\lfloor z \right\rfloor$ and $y = \text{frac}(z)$, respectively. Then, using equation~\eqref{eq_14} we have
\[
\arctan(\lfloor z \rfloor  + \text{frac}(z)) = \arctan\lfloor z \rfloor + \arctan \frac{\text{frac}(z)}{1 + (\lfloor z \rfloor + \text{frac}(z))\lfloor z \rfloor}
\]
or
\[
\arctan z = \arctan\lfloor z \rfloor + \arctan\frac{z - \lfloor z \rfloor}{1 + z\lfloor z \rfloor}
\]
from which it follows that
\[
\arctan\frac{1}{z}= 
\begin{cases}
\frac{\pi}{2} - \arctan\frac{z - \lfloor z \rfloor}{1 + z\lfloor z \rfloor}, & \text{if\; $z \in (0,1)$}; \\
{\text{undefined}}, & \text{if\; $z = 0$;} \\
\arctan\frac{1}{\lfloor z \rfloor} - \arctan\frac{z - \lfloor z \rfloor}{1 + z\lfloor z \rfloor}, & \text{otherwise,}
\end{cases}
\]
since
$$
\arctan z =
\begin{cases}
\frac{\pi}{2} - \arctan\frac{1}{z}, & \text{if\; $z > 0$;}\\
0, & \text{if\; $z = 0$;} \\
-\frac{\pi}{2} - \arctan\frac{1}{z}, & \text{otherwise.}
\end{cases}
$$
Consequently, we can write
\[
\arctan\frac{1}{z} = \arctan\frac{1}{\lfloor z \rfloor} - \arctan\frac{z - \lfloor z \rfloor}{1 + z\lfloor z \rfloor}, \qquad z\notin [0,1)
\]
or
\begin{equation}\label{eq_15}
\arctan\frac{1}{z} = \arctan\frac{1}{\lfloor z \rfloor} + \arctan\frac{\lfloor z \rfloor - z }{1 + z\lfloor z \rfloor}, \qquad z \notin [0,1).
\end{equation}

\begin{figure}[!ht]
\begin{center}
\includegraphics[width=24pc]{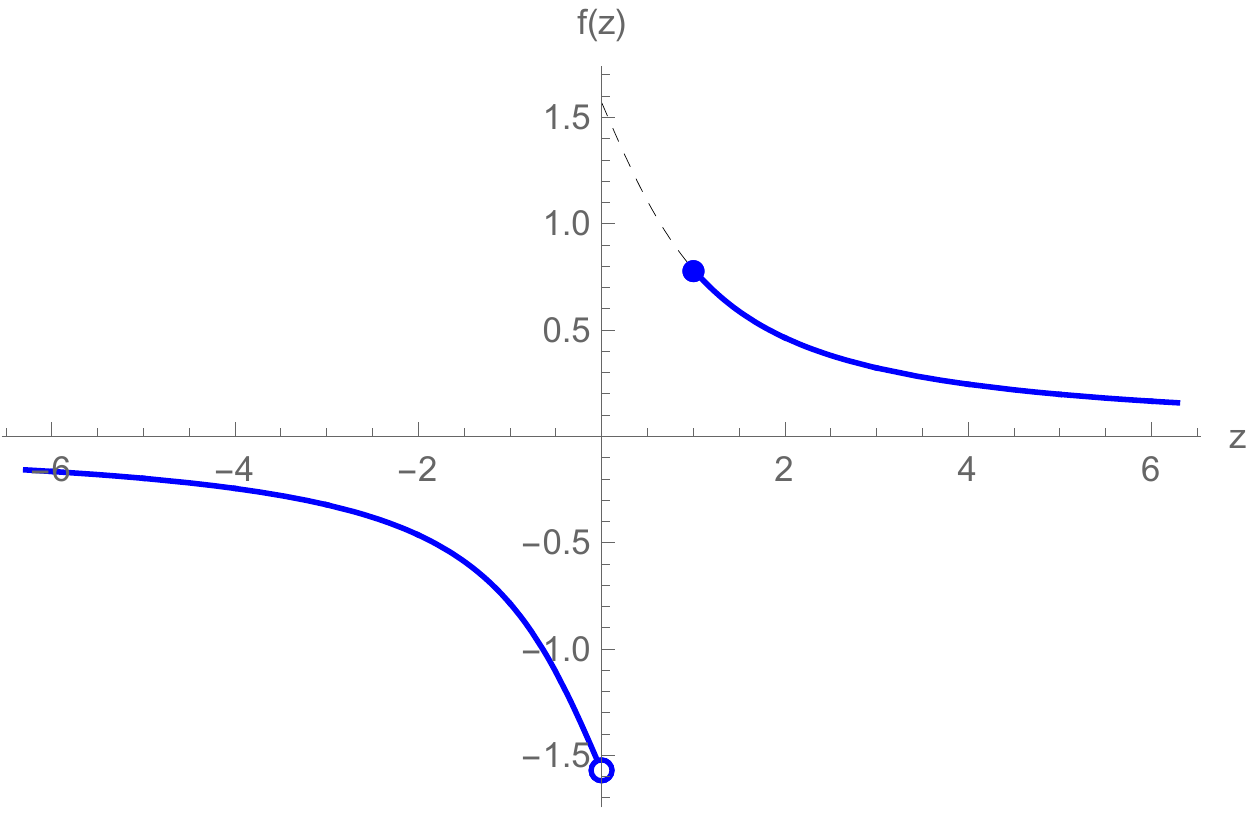}\hspace{2pc}
\caption{Function $f(z)$ undefined at $0\le z < 1$.}\label{fig1}
\end{center}
\end{figure}

Fig.~\ref{fig1} shows the function
$$
f(z) = \arctan\frac{1}{\lfloor z \rfloor} + \arctan \frac{\lfloor z \rfloor - z}{1 + z\lfloor z \rfloor},
$$
defined according to the right side of equation~\eqref{eq_15}. Open and filled circles along broken blue curve  in Fig.~\ref{fig1} indicate two points $(0,-\pi/2)$ and $(1,\pi/4)$, where the function $f(z)$ interrupts and resumes again. The dashed black curve, where the function $f(z)$ is not valid due to restriction $z\notin [0,1)$, is also shown for clarity.

Equation~\eqref{eq_15} is essential for generating a new form of the Machin-like formula for $\pi$. Consider a few examples showing how we can apply equation~\eqref{eq_15}. Taking $k = 6$ and using equation~\eqref{eq_9} we can find that
$$
\beta_1 = \left\lfloor\frac{\sqrt{2 + \sqrt{2 + \sqrt{2 + \sqrt{2 + \sqrt{2 + \sqrt 2}}}}}}{\sqrt{2 -\sqrt{2 + \sqrt{2 + \sqrt{2 + \sqrt{2 + \sqrt 2}}}}}}\right\rfloor = 40.
$$
With $\beta_1 = 40$ the equation~\eqref{eq_13} based on two-step iteration~\eqref{eq_12} leads to a quotient
\[
\beta_2 = -\frac{2634699316100146880926635665506082395762836079845121}{38035138859000075702655846657186322249216830232319}
\]
consisting of $52$ and $50$ digits in its numerator and denominator, respectively. Consequently, according to equation~\eqref{eq_7}, we can generate the two-term Machin-like formula for $\pi$
\begin{equation}\label{eq_16}
\frac{\pi}{4} = 32\arctan\frac{1}{40} + \arctan\frac{1}{\beta_2}.
\end{equation}

The two-term Machin-like formula~\eqref{eq_16} for $\pi$ has the Lehmer measure $\mu \approx 1.16751$. Although this value is smaller than the Lehmer measure $\mu \approx 1.51244$ corresponding to Chien-Lih's equation above, its application for the computation of $\pi$ may not be more efficient due to quotient consisting of many digits in its numerator and denominator. As mentioned already, this problem occurs due to exponentiation with increasing orders that is needed for computation of the arctangent function by using a series expansion like~\eqref{eq_6}. Therefore, it would be very desirable to obtain a Machin-like formula for $\pi$ with a reduced number of the digits in the numerator and an increased number of the digits in the denominator. Application of equation~\eqref{eq_15} may be one of the efficient ways to resolve such a problem.

Assuming that $z = \beta_2$ and using equation~\eqref{eq_15}, we can rewrite equation~\eqref{eq_16} as
\begin{equation}\label{eq_17}
\frac{\pi}{4} = 32\arctan\frac{1}{40} - \arctan\frac{1}{70} + \arctan\frac{1}{\beta_3},
\end{equation}
where the quotient
$$
\beta_3 = -\frac{184466987265869281740567152432082954025647742419390789}{27760404029858418259273600496960161682342036417209},
$$
consists of $54$ and $50$ digits in its numerator and denominator. The three-term Machin-like formula~\eqref{eq_17} for $\pi$ does not look interesting since its second integer $70$ has the same order as the first integer $40$. However, if we repeat again the same procedure for equation~\eqref{eq_17} by assuming that $z = \beta_3$, then using equation~\eqref{eq_15} we can observe that in the next 
Machin-like formula for $\pi$
\begin{equation}\label{eq_18}
\frac{\pi}{4} = 32\arctan\frac{1}{40} - \arctan\frac{1}{70} - \arctan\frac{1}{6645} + \arctan\frac{1}{\beta_4},
\end{equation}
where
$$
\beta_4 = -\frac{612891579071052703512243493592395863230295465359444105057}{448756269953796152961435108660176757544786481508},
$$
the third constant $6645$ is by two order of the magnitude larger than the second integer $70$ (the quotient $\beta_4$ consists of $57$ and $48$ digits in its numerator and denominator). Repeating same procedure again for equation~\eqref{eq_18} by assuming that $z = \beta_4$,
with help of equation~\eqref{eq_15} we get
\begin{equation}\label{eq_19}
\begin{aligned}
\frac{\pi}{4} &= 32\arctan\frac{1}{40} - \arctan\frac{1}{70} - \arctan\frac{1}{6645} \\
&- \arctan\frac{1}{1365756025} + \arctan\frac{1}\beta_5,
\end{aligned}
\end{equation}
where the quotient
$$
\beta_5 = \scalebox{1.2}[1.7]{$-\frac{837060366788054133363141482594659697353287103005016334677117199933}{374870864016658098706770220951460879098657980643}$}
$$
consists of $66$ and $48$ digits in its numerator and denominator, respectively.

As we can see, the fourth integer in equation~\eqref{eq_19}, $1365756025$, is significantly larger than the third integer $6645$. Repeating the same procedure over and over again, we notice that each next integer is larger than the previous one by many orders of magnitude.

Let us show how this methodology can also be used to modify the Machin-like formulas for $\pi$. Consider, as an example, the following identity that was found by Wetherfield in 2004~\cite{Wetherfield2004}
\[
\begin{aligned}
\frac{\pi}{4} &= 83\arctan\frac{1}{107} + 17\arctan\frac{1}{1710} - 22 \arctan\frac{1}{103697} \\
&-12\arctan\frac{2}{2513489} - 22\arctan\frac{2}{18280007883}.
\end{aligned}
\]
Although the Lehmer measure of this formula is small $\mu \approx 1.26579$, the arguments in the last two arctangent functions are not the integer reciprocals. However, applying equation~\eqref{eq_15} in a sequence first for
$$
z=-\frac{2513489}{2}
$$
and then for
$$
z=-\frac{18280007883}{2},
$$
we obtain the identity
\[
\begin{aligned}
\frac{\pi}{4} &= 83\arctan\frac{1}{107} + 17\arctan\frac{1}{1710} - 22 \arctan\frac{1}{103697} \\
&-12\arctan\frac{1}{1256744} - 22\arctan\frac{1}{9140003941} \\
&+12\arctan\frac{1}{3158812219818} + 22\arctan\frac{1}{167079344092131066905},
\end{aligned}
\]
where all arguments are the integer reciprocals. The Lehmer measure for this identity is $\mu \approx 1.39524$.

\subsection{Generalization}

It is not difficult to see by induction that the procedure described above can be generalized as a new form of the Machin-like formula for $\pi$
\begin{equation}\label{eq_20}
\frac{\pi}{4} = 2^{k - 1}\arctan\frac{1}{A_k} + \left(\sum\limits_{m = 1}^M \arctan\frac{1}{\lfloor{B_{m,k}}\rfloor}\right) + \arctan\frac{1}{B_{M + 1,k}},
\end{equation}
where $A_k$ is the integer defined by equation~\eqref{eq_9} and
\begin{equation}\label{eq_21}
B_{m,k} = \frac{1 + \lfloor{B_{m - 1,k}}\rfloor {B_{m - 1,k}}}{\lfloor{B_{m - 1,k}}\rfloor - B_{m - 1,k}}, \qquad m \ge 2
\end{equation}
with initial integer $B_{1,k}$ that can be computed either by using equation~\eqref{eq_11} or, more efficiently, by using equation~\eqref{eq_13} based on the two-step iteration~\eqref{eq_12}.

Algorithmic implementation of equation~\eqref{eq_20} implies two important rules. First, since the integer $B_{0,k}$ is not defined, it follows that at $M = 0$ the sum of arctangent functions
$$
\left.\sum\limits_{m = 1}^M \arctan\frac{1}{\lfloor B_{m,k} \rfloor}\right|_{M = 0} = 0.
$$
Second, if
$$
\lfloor B_{M + 1,k} \rfloor - B_{M + 1,k} = 0,
$$
then no further iteration is required as the fractional part of the number $B_{M + 1,k}$ does not exist.

As a simplest example, consider $k = 4$. Using equation~\eqref{eq_9} we can find that
$$
A_4 = \left\lfloor \frac{\sqrt {2 + \sqrt{2 + \sqrt{2 + \sqrt 2}}}}{\sqrt {2 - \sqrt{2 + \sqrt{2 + \sqrt 2}}}}\right\rfloor = 10.
$$
Then, using two-step iteration~\eqref{eq_13}, we have
\begin{align*}
\begin{aligned}
\sigma_1 &= \frac{A_4^2 - 1}{A_4^2 + 1} = \frac{99}{101}, &\tau_1 &= \frac{2A_4}{A_4^2 + 1} = \frac{20}{101}, \\
\sigma_2 &= \sigma_1^2 - \tau_1^2 = \frac{9401}{10201}, &\tau_2 &= 2\sigma_1 \tau_1 = \frac{3960}{10201}, \\
\sigma_3 &= \sigma_2^2 - \tau_2^2 = \frac{72697201}{104060401}, &\tau_3 &= \;2\sigma_2 \tau_2 = \frac{74455920}{104060401}, \\
\sigma_4 &= \sigma_3^2 - \tau_3^2 =  -\frac{258800989811999}{10828567056280801}, &\tau_4 &= 2\sigma_3 \tau_3 = \frac{10825473963759840}{10828567056280801}.
\end{aligned}
\end{align*}

Substituting $\sigma_4$ and $\tau_4$ into equation~\eqref{eq_14} leads to
$$
B_{1,4} = \frac{\sigma _4}{{1 - \tau_4}} = -\frac{147153121}{1758719}.
$$
Therefore, equation~\eqref{eq_8} yields
\begin{equation}\label{eq_22}
\frac{\pi}{4} = 8\arctan\frac{1}{10} - \arctan\frac{1758719}{147153121}.
\end{equation}
As we can see,  at $M = 0$ the identity~\eqref{eq_22} is consistent with equation~\eqref{eq_20}.

Using $k = 4$ and $M = 2$ in equation~\eqref{eq_20}, we obtain
\begin{equation}\label{eq_23}
\begin{aligned}
\frac{\pi}{4} &= 8\arctan\frac{1}{10} - \arctan\frac{1}{84} - \arctan\frac{1}{21342}\ \\
&- \arctan\frac{266167}{263843055464261}.
\end{aligned}
\end{equation}
The argument of the last arctangent function in equation~\eqref{eq_23} is not an integer reciprocal. Therefore, we can use iteration formula~\eqref{eq_21} again. However, at $k = 4$ and $M = 5$, the equation~\eqref{eq_20} results in
\begin{align}\label{eq_24}
\frac{\pi}{4} &= 8\arctan\frac{1}{10} - \arctan\frac{1}{84} - \arctan\frac{1}{21342} \nonumber \\
&- \arctan\frac{1}{991268848} - \arctan\frac{1}{193018008592515208050} \nonumber \\
&- \arctan\frac{1}{197967899896401851763240424238758988350338} \\
&- \arctan\scalebox{0.9}[1.5]{$\frac{1}{117573868168175352930277752844194126767991915008537018836932014293678271636885792397}$}\vspace{0.25cm}, \nonumber
\end{align}
where the largest integer consists of $84$ digits. Consequently, no further iteration is required since the last argument of the arctangent function is an integer reciprocal now. As we can see, starting from the second arctangent term, the integers in equation~\eqref{eq_24} increase by many orders of the magnitude at each step of the iteration.

To generate the multi-term Machin-like formulas for $\pi$ with only the integer reciprocals, a method known as Todd's process is commonly applied~\cite{Wetherfield1996,Arndt2011,Todd1949}. Generally, Todd's process is quite complicated and requires matrix manipulations based on a set of  primes. However, considering equation~\eqref{eq_24} as an example, we can conclude that the proposed iterative method can be used as a simple alternative to Todd's process to generate the multi-term Machin-like formulas for $\pi$ consisting of only integer reciprocals.

\subsection{Approximation}

Since iteration~\eqref{eq_21} leads to 
\begin{equation}\label{eq_25}
|\lfloor B_{m,k} \rfloor|  \gg |\lfloor B_{m - 1,k} \rfloor| \gg |\lfloor B_{m - 2,k} \rfloor| \gg \cdots \gg |\lfloor B_{2,k} \rfloor|,
\end{equation}
we can infer that ${B_{M + 1,k}}$ is the largest by absolute value. Consequently, it is reasonable to approximate the last arctangent function term in equation~\eqref{eq_20} as
\begin{equation}\label{eq_26}
\arctan\frac{1}{B_{M + 1,k}} \approx \frac{1}{B_{M + 1,k}}, \qquad \text{if}\;B_{M + 1,k}\ne \lfloor B_{M + 1,k} \rfloor
\end{equation}
in accordance with the Maclaurin series expansion~\eqref{eq_6}. Thus, we can write the following approximation
\begin{equation}\label{eq_27}
\frac{\pi}{4} \approx 2^{k - 1}\arctan\frac{1}{A_k} + \left(\sum\limits_{m = 1}^M \arctan\frac{1}{\lfloor B_{m,k} \rfloor}\right) + \frac{1}{B_{M + 1,k}}.
\end{equation}

Table~\ref{tab1} shows quantity of correct digits depending on the integer $M$ in equation~\eqref{eq_27}. As we can see from this table, starting from $M = 2$, each increment of the integer $M$ in approximation~\eqref{eq_27} doubles the number of correct digits of $\pi$.  Consequently, we can estimate that, at $M = 26$, the approximation~\eqref{eq_27} can provide more than a billion digits of $\pi$.
\begin{table}[!ht]
\begin{center}
\begin{tabular}{c|c}
Integer $M$ & Correct digits of $\pi$ \\
\hline
0  & 5     \\
1  & 11    \\
2  & 27    \\
3  & 54    \\
4  & 98    \\
5  & 222   \\
6  & 444   \\
7  & 889   \\
8  & 1783  \\
9  & 3567  \\
10 & 7136  \\
11 & 14273 \\
12 & 28546
\end{tabular}
\caption{Correct digits of $\pi$ for different $M$ in equation~\eqref{eq_27}.}\label{tab1}
\end{center}
\end{table}

Since the initial values $A_k$ and $B_{1,k}$ are larger with increasing $k$, we may reduce the number of the summation terms in approximation~\eqref{eq_27}. For example, even at relatively small value $k = 17$, a computational test shows that at $M = 0$ in equation~\eqref{eq_27}, the number of correct digits of $\pi$ is $19$. Therefore, by doubling digits after each step of iteration, we can estimate that at $k = 17$ and $M = 24$, the approximation~\eqref{eq_27} can provide more than a billion digits of $\pi$.

Application of the approximation~\eqref{eq_27} may be advantageous. Due to the large absolute magnitude of the number $B_{M+1,k}$, the last arctangent function can be replaced by its argument in accordance with equation~\eqref{eq_27}. As a result, we do not need to include the quotient $B_{M+1,k}$ in computation of the Lehmer measure~\eqref{eq_7}. Moreover, this approach does not require exponentiation in the computation of the arctangent function with a problematic quotient and, if we need to improve the accuracy, we can increase the integer $M$ in approximation~\eqref{eq_27}.

It should be noted that increasing the number of the summation terms in equation~\eqref{eq_27} increases the Lehmer measure $\mu$. However, our empirical results show that at $k\ge 17$, only the first $18$ terms actually contribute for $\mu$. Due to condition~\eqref{eq_25} the contribution of the additional arctangent function terms to the Lehmer measure of equation~\eqref{eq_27} becomes vanishingly small. For example, at $k = 17$ and $M = 0$, equation~\eqref{eq_27} gives Lehmer's measure $\mu  \approx 0.203195$. The value $\mu$ increases with increasing $M$ since more and more arctangent function terms are added as $M$ increases. However, after $18$  steps of iteration, the Lehmer measure reaches the value $\mu  \approx 0.50222$ and further remains practically unchanged with increasing $M$. This stabilization is due to rapidly increasing magnitude of integers $\left\lfloor B_{m,k} \right\rfloor$ that becomes particularly evident at $k\ge 17$. Therefore, the presence of very large integers $\left\lfloor B_{m,k} \right\rfloor$ and exclusion of the arctangent function $\arctan\left(1/B_{M+1,k}\right)$ indicate that the approximation~\eqref{eq_27} may be promising for computing $\pi$ with rapid convergence.

\subsection{Arctangent function}

Apart from the Maclaurin series expansion~\eqref{eq_6}, the following limit~\cite{Abrarov2021}
\[
\arctan x = \lim_{N\to\infty}\sum_{n = 1}^N \frac{N x}{N^2 + (n - 1)nx^2}
\]
can also be used for computation of the arctangent function terms. Although both these equations are simple, their implementation provides relatively slow convergence. There are several other interesting equations for the arctangent function that do not need irrational (surd) numbers in computation~\cite{Olds1963,Henrici1977,Milgram2006,Lorentzen2008,Sofo2012}. However, the following two series expansions~\cite{Castellanos1988, Chien-Lih2005}
\begin{equation}\label{eq_28}
\arctan x = \frac{x}{1+x^2}\;_2F_1\left(1,1;\frac{3}{2};\frac{x^2}{1+x^2}\right) = \sum\limits_{n = 0}^\infty\frac{2^{2n}(n!)^2}{(2n + 1)!}\frac{x^{2n + 1}}{(1 + {x^2})^{n + 1}},
\end{equation}
where $_2F_1(a,b;c;z)$ denotes the hypergeometric function, and
\begin{equation}\label{eq_29}
\arctan x = 2\sum\limits_{n = 1}^\infty\frac{1}{2n - 1}\frac{g_n(x)}{g_n^2(x) + h_n^2(x)},
\end{equation}
where the expansion coefficients are computed by iteration
$$
g_1(x) = 2/x,\,\,\,h_1(x) = 1,
$$
$$
g_n(x) = g_{n - 1}(x)(1 - 4/x^2) + 4h_{n - 1}(x)/x,
$$
$$
h_n(x) = h_{n - 1}(x)(1 - 4/x^2) - 4g_{n - 1}(x)/x,
$$
are found to be most suitable for computation due to their rapid convergence.

Chien-Lih showed a simple and elegant derivation of the Euler's series expansion~\eqref{eq_28} by taking the integral
$$
\arctan x = \int_0^{\pi/2}\frac{x\sin u}{1 + x^2}\;\frac{1}{\left(1 - \frac{x^2\sin^2 u}{1 + x^2}\right)}du
$$
in terms of geometric series~\cite{Chien-Lih2005}
$$
\frac{1}{1 - \frac{x^2\sin^2 u}{1 + x^2}}=\sum_{n = 0}^{\infty}\frac{x^{2n}\sin^{2n}u}{(1 + x^2)^n}.
$$
The series expansion~\eqref{eq_29} represents a trivial rearrangement of the equation
\[
\arctan x = i\sum_{n=1}^{\infty}\frac{1}{2n - 1}\left(\frac{1}{(1 + 2i/x)^{2n - 1}} - \frac{1}{(1 - 2i/x)^{2n - 1}}\right)
\]
that was derived in~\cite{Abrarov2018}. The computational test we performed shows that equation~\eqref{eq_29} is more rapid in convergence than equation~\eqref{eq_28}. Therefore, the application of the iteration-based series expansion~\eqref{eq_29} may be more preferable for computation of the arctangent function terms in the approximation~\eqref{eq_27}~\cite{Abrarov2017}.

\section{Alternative method}

It is convenient to use calligraphic letters $\mathcal{A}$ and $\mathcal{B}$ to keep consistency with equation~\eqref{eq_20}. Using the methodology described above in the section~\ref{RD}, we can also derive the Machin-like formula for $\pi$ in an alternative form as
\begin{equation}\label{eq_30}
\frac{\pi}{4} = 2^{k - 1}\left(\left(\sum\limits_{m = 1}^M\arctan\frac{1}{\lfloor{\mathcal A}_{m,\ell ,k}\rfloor}\right) + \arctan\frac{1}{{\mathcal A}_{M + 1,\ell ,k}}\right) + \arctan\frac{1}{{\mathcal B}_{\ell,k}},
\end{equation}
where the constants ${{\mathcal A}_{m,\ell ,k}}$ can be computed by iteration
\begin{equation}\label{eq_31}
{\mathcal A}_{m,\ell ,k} = \frac{1 + {\mathcal A}_{m - 1,\ell,k}\lfloor{\mathcal A}_{m - 1,\ell ,k}\rfloor}{\lfloor{\mathcal A}_{m - 1,\ell ,k}\rfloor - {\mathcal A}_{m - 1,\ell ,k}}
\end{equation}
with initial number defined as
\begin{equation}\label{eq_32}
{\mathcal A}_{1,\ell,k} = \frac{1}{10^\ell}\left\lfloor 10^\ell\frac{a_k}{\sqrt{2 - a_{k - 1}}}\right\rfloor.
\end{equation}

Similar to the Machin-like formula~\eqref{eq_20} for $\pi$ this equation also implies the same two rules. Since ${\mathcal A}_{k,\ell,0}$ is not defined, we imply that
\[
\left.\sum\limits_{m = 1}^M \arctan\frac{1}{\lfloor{\mathcal A}_{m,\ell,k}\rfloor}\right|_{M = 0} = 0
\]
and if the following condition
$$
{\mathcal A}_{M + 1,\ell,k} - \lfloor{\mathcal A}_{M + 1,\ell,k}\rfloor = 0
$$
is satisfied, further iteration is not needed since the number $\mathcal A_{M + 1,\ell,k}$ is an integer.

Consider the following examples. At $\ell=2$ and $k=4$, from equation~\eqref{eq_32} it follows that the initial number is
$$
{\mathcal A}_{1,2,4} = \frac{1}{10^2}\left\lfloor 10^2\frac{a_4}{\sqrt{2 - {a_3}}}\right\rfloor = \frac{1}{10^2}\left\lfloor 10^2\frac{\sqrt {2 + \sqrt{2 + \sqrt{2 + \sqrt 2}}}}{\sqrt{2 - \sqrt{2 + \sqrt{2 + \sqrt 2}}}}\right\rfloor = \frac{203}{20}.
$$
Consequently, substituting ${\mathcal A}_{1,2,4}$ into equation~\eqref{eq_11} for $\beta_1$ or using equation~\eqref{eq_13} based on two-step iteration~\eqref{eq_12} we can find that
$$
{\mathcal B}_{\ell,k} = -\frac{4239006656613482881}{1033248635280959}.
$$

Thus, at $M=0$, $\ell=2$ and $k=4$ equation~\eqref{eq_30} leads to
\begin{equation}\label{eq_33}
\frac{\pi}{4} = 8\arctan\frac{20}{203} - \arctan\frac{1033248635280959}{4239006656613482881}.
\end{equation}
At $M=3$, $\ell=2$ and $k=4$ equation~\eqref{eq_30} gives
\begin{equation}\label{eq_34}
\begin{aligned}
\frac{\pi}{4} &= 8\left(\arctan\frac{1}{10} - \arctan\frac{1}{684} - \arctan\frac{2}{1402203}\right) \\
&-\arctan\frac{1033248635280959}{4239006656613482881}.
\end{aligned}
\end{equation}
Since in equation~\eqref{eq_34}
$$
{\mathcal A}_{3,2,4} = -\frac{1402203}{2}
$$
is not an integer, we can apply iteration formula~\eqref{eq_31} again. This leads to
\begin{equation}\label{eq_35}
\begin{aligned}
\frac{\pi}{4} &= 8\left(\arctan\frac{1}{10} - \arctan\frac{1}{684} - \arctan\frac{1}{701102}\right.\\
&\left. -\arctan\frac{1}{983087327708}\right) - \arctan\frac{1033248635280959}{4239006656613482881}.
\end{aligned}
\end{equation}
As we can see
$$
\mathcal{A}_{4,2,4} = -983087327708
$$
is an integer now. Therefore, no further iteration is required.

The last term in equation~\eqref{eq_35}
\[
\begin{aligned}
\arctan\frac{1}{\mathcal{B}_{2,4}} &= -\arctan\frac{1033248635280959}{4239006656613482881} \\
&=\arctan \left(-\frac{1033248635280959}{4239006656613482881}\right)
\end{aligned}
\]
can also be represented as a sum of the arctangent functions with the integer reciprocals by using the identity~\eqref{eq_15}.

When $\ell$ or $k$ increases, the constant ${\mathcal B}_{\ell,k}$ also increases by absolute value. Therefore, according to series expansion~\eqref{eq_6}, we can write
$$
\arctan\frac{1}{{\mathcal B}_{\ell,k}} \approx \frac{1}{{\mathcal B}_{\ell,k}}, \qquad \ell \gg 1 \;\text{or}\; k \gg 1
$$
and modify equation~\eqref{eq_30} as
\begin{equation}\label{eq_36}
\frac{\pi}{4} \approx 2^{k - 1}\left(\left(\sum\limits_{m = 1}^M \arctan\frac{1}{\lfloor {\mathcal A}_{m,\ell,k}\rfloor}\right) + \arctan\frac{1}{{\mathcal A}_{M + 1,\ell,k}}\right) + \frac{1}{{\mathcal B}_{\ell,k}}
\end{equation}
to approximate $\pi$. However, in contrast to equation~\eqref{eq_27} approximation~\eqref{eq_36} does not improve accuracy with increasing $M$ since the constant ${\mathcal B}_{\ell,k}$ is independent of $M$.

\section{Supplement}

The file \href{https://cs.uwaterloo.ca/journals/JIS/VOL25/Abrarov/supplement.txt}{\textcolor{blue}{\it supplement.txt}} provides Mathematica programs that can be copy-pasted to the Mathematica notebook to validate the main results obtained in this study.

\section{Conclusion}

In this work, we propose a new form of the Machin-like formula~\eqref{eq_20} for $\pi$ that is generated by using iteration formula~\eqref{eq_21}. Due to condition~\eqref{eq_25}, the application of this form of the Machin-like formula may be promising for computation of the constant $\pi$ with rapid convergence. Approximation~\eqref{eq_27} shows that, at $k\ge 17$, the Lehmer measure remains small and practically does not increase after $18$ steps of iteration.

\section*{Acknowledgments}

This work is supported by National Research Council Canada, Thoth Technology Inc., York University and Epic College of Technology.


\end{document}